\documentclass{article}

\usepackage[unicode=true,pdfusetitle,
bookmarks=true,bookmarksnumbered=true,bookmarksopen=true,bookmarksopenlevel=2,
breaklinks=false,pdfborder={0 0
  1},backref=false,colorlinks=false]{hyperref}
\usepackage{amsfonts, amsmath, amssymb}
\usepackage{graphicx}
\usepackage{lmodern}
\usepackage{verbatim}
\usepackage{multirow}
\usepackage{natbib}
\usepackage{booktabs}
\usepackage{setspace}
\usepackage[margin=1in]{geometry}

\newcommand\BibTeX{{\rmfamily B\kern-.05em \textsc{i\kern-.025em b}\kern-.08em
T\kern-.1667em\lower.7ex\hbox{E}\kern-.125emX}}

\renewcommand{\c}[1]{\texttt{#1}}
\renewcommand{\b}[1]{\textbf{#1}}

\newcommand{\R}{\mathbb{R}}
\newcommand{\M}{{\cal M}}

\newcommand{\bX}{{\bf X}}

\newcommand{\bx}{{\bf x}}
\newcommand{\by}{{\bf y}}
\newcommand{\bz}{{\bf z}}
\newcommand{\bu}{{\bf u}}

\newcommand{\bbf}{{\bf f}}
\newcommand{\tbx}{{\tilde{\bx}}}

\newcommand{\amax}[1]{\underset{#1}{\textup{argmax} \,}}

\title{\textbf{Combining clustering of variables and feature selection
    using random forests}}
\author{Marie {\sc Chavent}$^{1,2}$, Robin {\sc Genuer}$^{3, 4}$,
  J\'er\^ome {\sc Saracco}$^{1,2}$}
\date{}

\begin{document}

\maketitle

\begin{center}{\small
$^1$ Institut de Math\'ematiques de Bordeaux, UMR CNRS 5251,
Universit\'e de Bordeaux, \\
351 cours de la lib\'eration,   33405 Talence Cedex, France. \\

\smallskip

$^2$ INRIA Bordeaux Sud-Ouest, CQFD team, France.\\

\smallskip

$^3$ ISPED, INSERM U-1219, Universit\'e de Bordeaux,\\
146 rue L\'eo Saignat, 33076 Bordeaux, France.\\

\smallskip

$^4$ INRIA Bordeaux Sud-Ouest, SISTM team, France.\\
}

\smallskip

{\small e-mail~: {\tt Marie.Chavent,Jerome.Saracco@math.u-bordeaux.fr}, {\tt Robin.Genuer@u-bordeaux.fr}}

\end{center}

\begin{abstract}
Standard approaches to tackle high-dimensional supervised classification problem often include variable selection and dimension reduction procedures. The novel methodology proposed in this paper combines clustering of variables and feature selection. More precisely, hierarchical clustering of variables procedure allows to build groups of correlated variables in order to reduce the redundancy of information and summarizes each group by a synthetic numerical variable. Originality is that the groups of variables (and the number of groups) are unknown a priori. Moreover the clustering approach used can deal with both numerical and categorical variables (i.e. mixed dataset). Among all the possible partitions resulting from dendrogram cuts, the most relevant synthetic variables (i.e. groups of variables) are selected with a variable selection procedure using random forests. Numerical performances of the proposed approach are compared with direct applications of random forests and variable selection using random forests on the original $p$ variables. Improvements obtained with  the proposed methodology are illustrated on two simulated mixed datasets (cases $n>p$ and $n<p$, where $n$ is the sample size) and on a real proteomic dataset. Via the selection of groups of variables (based on the synthetic variables), interpretability of the results becomes easier.
\end{abstract}



\section{Introduction}

This paper addresses the problems of dimension reduction and variable selection
in the context of supervised classification. In this framework, there are
often two objectives : prediction (to be able to predict classes
associated to new observations) and feature selection (to
extract the most interesting variables).
Typically in a medical context, the first goal could be to succeed
in predicting if
patients will or will not respond well to a treatment, given e.g. their gene
expression profile, whereas the second aim could be to determine which
part of the genome is responsible for the good or bad response to the
treatment. These two objectives can be related, because it may be easier to
perform prediction if useless variables have already been eliminated.
Those goals are indeed especially relevant when dealing with e.g. genomics or
proteomics data, where the number $p$ of variables largely exceeds the number
$n$ of available observations.

A classical way of addressing such issues is to use a variable
selection technique \citep[see e.g.][]{Guyon03, tibshirani1996regression,
zou2005regularization}. The hope is that the method will
be able to select the most interesting variables, while preserving good
prediction performances. Moreover, data often come
with many highly correlated variables and succeeding in selecting all
variables in a group of correlated variables can be very difficult.
Even if obtaining all these variables does not always appear to be interesting
for prediction purpose, it can be useful for interpretation purpose,
depending on the application.

We stress that the framework considered in this work is the \emph{supervised}
classification framework, in which an output variable is observed. Hence, we
do not address the problem of variable selection in a clustering (or
unsupervised classification) context as e.g. in \citet{raftery2006variable,
tadesse2005bayesian, maugis2009variable}.

In this article, we propose and evaluate a novel methodology for dimension
reduction and variable selection, which combines clustering of variables and
feature selection using random forests.
The clustering of variables method, introduced in \citep{Chavent12b} and
denoted by \b{CoV} hereafter, allows to eliminate redundancy of explanatory
variables. This clustering approach can deal with both numerical and
categorical variables (i.e. mixed dataset).
Note that this method clusters the \emph{variables}, which differs
from the (more classical) problem of \emph{individuals} clustering
\citep{kaufman2009finding, gordon1999classification}.
The clustering of variables groups together highly correlated variables and
provides for each group (cluster) a synthetic variable
which is a numerical variable summarizing the variables within the cluster.
The main advantage of this approach is to
 eliminate redundancy and to keep all the variables together in a cluster  during
the rest of the analysis. Moreover it reduces the dimension of the data
by replacing the $p$ original variables by $K$ synthetic variables (where $K$
denotes the selected number of clusters). In the proposed methodology the number $K$ of
synthetic variables is optimized according to prediction performance of a 
random forest classifier trained with those synthetic variables.
Let $K^*$ denote the optimal number of clusters of variables.
Note that the reduction of dimension provides $K^*$ synthetic variables which  only use the original variables within the selected clusters, unlike the principal components in principal component analysis (PCA) \citep{jolliffe2002principal}. Hence, in \b{CoV}, an original variable takes action in the construction of a unique synthetic variable, which makes the interpretation easier.
Once the dimension reduction is done, the most important synthetic variables
are then selected using a procedure based on random forests (\b{RF}), introduced in
\citet{Genuer10}. This variable selection
procedure, denoted \b{VSURF} hereafter, is applied to the reduced dataset
consisting of the $n$ observations described with the $K^*$ synthetic variables,
and leads to provide $m \leq K^*$ relevant synthetic variables.
Thus the prediction for new observations can be done
with a predictor built on these $m$ selected synthetic variables (i.e. a list of clusters of
variables).
Hence, the output of the proposed method is a classifier and a set of selected variables
with the additional information of the group structure of those variables.
This additional information can be of great interest depending of the application,
and as far as we know, this type of methodology has not already been
introduced in the literature.

Note that the proposed approach \emph{does not} require
definition of a priori groups of variables. Hence, it is different from
(sparse) group lasso or (sparse) group partial least squares (PLS) approaches
\citep{yuan2006model, simon2013sparse, liquet2015group}.
In addition, our method also differs from sparse PLS techniques
\citep{chun2010sparse, le2011sparse}, which perform dimension reduction and
variable selection without any group structure information.

In the proposed methodology combining \b{CoV} and \b{VSURF}, the synthetic variables are linear combinations of the variables within a group, and the \b{VSURF}
procedure is purely non parametric. This avoids to make assumptions
between the response variable and synthetic variables (which are unknown at the
beginning of the algorithm).
Moreover, since we heavily use the classifier (to choose
the optimal number $K^*$ of synthetic variables and to select the $m$ most important of them),
\b{RF} is chosen as an easy-to-use and well-performing algorithm
\citep{fernandez2014we}.

Finally, let us stress that the main objective of this paper is to introduce
a novel procedure which 
\begin{itemize}
\item[-] selects groups of
informative variables, 
\item[-] can deal with mixed data, 
\item[-] provides good numerical
performance, 
\item[-] and eases interpretation of results.
\end{itemize}

\medskip

The rest of the paper is organized as follows.
Sections~\ref{subsec:clustofvar} and \ref{subsec:vsurf} give an overview of both the clustering
of variables method and the feature selection procedure. 
The methodology combining \b{CoV} and \b{VSURF} is then described in  Section~\ref{subsec:covsurf}.
In Section~\ref{sec:simu-study}, the numerical performances of this
methodology are compared, via a simulation study, with a straightforward applications of \b{VSURF} or \b{RF}
(on the original $p$ variables) and with application of \b{RF} on all the synthetic variables given by \b{CoV}.
Two cases are addressed:  $n>p$ and $n<p$.
A real proteomic
dataset is analyzed in Section~\ref{sec:real-data}.
Finally, let us mention that the proposed methodology is available and implemented in the R package {\tt CoVVSURF}\footnote{\url{https://github.com/robingenuer/CoVVSURF}}  including the simulation procedure exhibited in Section~\ref{sec:simu-study}. 

\section{Description of the methodology}
\label{sec:brief-descr}

First, the two underlying methods, \b{CoV} for unsupervised dimension reduction and \b{VSURF} for both variable selection and prediction, are presented  respectively in Sections~\ref{subsec:clustofvar} and \ref{subsec:vsurf}. 
Then, the proposed methodology combining \b{CoV} and \b{VSURF}, named \b{CoV/VSURF} hereafter, is  described in  Section~\ref{subsec:covsurf}.

Let us consider a $p$-dimensional explanatory variable $X=(X^1,\dots,X^j,\dots,X^p)'$ and a univariate response variable $Y$, which takes its values in $\{1, \ldots, L\}$.
Let  $n$ be the number of observations of these variables.
More precisely, let us consider a
set of  $p_1$ numerical variables measured on the $n$ observations denoted by $\{\bx^1,\dots,\bx^{p_1}\}$, and a set of  $p_2$ categorical variables denoted by $\{\tbx^1,\dots,\tbx^{p_2}\}$ with $p_1+p_2=p$.  Let $\bX = [\bx^1,\dots,\bx^{p_1}, \tbx^1,\dots,\tbx^{p_2}]$  be the corresponding  data matrix of dimension $n \times p$.
The $i$-th row of $\bX$ is denoted $\bx_i$.
Let $\by$ be the vector of the $n$ observations of the response variable.

\subsection{Clustering of variables}
\label{subsec:clustofvar}

The objective of clustering of variables is to sort  variables into homogeneous clusters, that is to construct clusters of variables which are strongly related to each other and thus provide similar information. The idea is to summarize all  the variables belonging to a cluster by a synthetic numerical variable which is the most ``linked'' to all the variables within this cluster.
In this section, we focus on a clustering method  
based on the  principal component analysis method \b{PCAmix} \citep{kiers1991simple, Chavent12a} defined for a mixture of categorical and numerical variables.
More precisely, we present the ascendant hierarchical clustering algorithm implemented in the R \citep{Rcoreteam} package {\tt ClustOfVar} \citep{Chavent12b}.

Note that like PCA, \b{CoV} is a dimension reduction method  but, contrary to PCA, it can be helpful for variable selection. Indeed, each synthetic variable of \b{CoV} is a linear combination of a subset of variables (the variables within the corresponding cluster) whereas the principal components in PCA are linear combination of all the original variables. Hence, selecting synthetic variables of \b{CoV} means selecting of subsets of original variables, which is not the case when selecting principal components in PCA.

\medskip

\noindent
\textbf{Synthetic variable of a cluster $C_k$}

\noindent
This variable is defined as the numerical variable $\bbf^k \in \R^n$  which is the ``most linked'' to all the variables in $C_k$:  
\begin{equation}
\bbf^k=\displaystyle \arg\max_{\bu \in \R^n} \left\{\sum_{\bx^j \in C_k} r^2_{\bu,\bx^j} + \sum_{\tbx^j \in C_k} \eta^2_{\bu|\tbx^j}  \right\},
\label{fkdef}
\end{equation}
where $r^2_{\bu,\bx^j}\in [0,1]$ is the squared Pearson correlation between the numerical variables $\bu$ and $\bx^j$,   
and $\eta^2_{\bu|\tbx^j}\in [0,1]$ is the correlation ratio between $\bu$ and the categorical variable $\tbx^j$ (which measures the part of the variance of $\bu$ explained by the levels of $\tbx^j$). 
It has been shown that:
\begin{itemize}
\item $\bbf^k$ is the first principal component of  the  \b{PCAmix} method applied to  the variables in $C_k$. 
\item  $\bbf^k$ is  a linear combination of the numerical variables  and of the dummy (indicator) variables of the  levels  of the categorical variables in $C_k$.
This linear combination can  be used to predict the value (score)  of a new observation on the  synthetic variable of $C_k$.
\end{itemize}
Details on the \b{PCAmix} method and on the prediction of principal component scores can be found in the Appendix. Note that the \b{PCAmix} algorithm is implemented in the R package  {\tt PCAmixdata} \citep{chavent2014multivariate}.

\medskip

\noindent
\textbf{Homogeneity $H$ of a cluster $C_k$}

\noindent
The following criterion $H$ measures the adequacy between the variables within the cluster  and its associated synthetic variable $\bbf^k$. It is defined as follows:
\begin{equation}
H(C_k)=\sum_{\bx^j \in C_k} r^2_{\bx^j,\bbf^k} + \sum_{\tbx^j \in C_k} \eta^2_{\bbf^k|\tbx^j}=\lambda_1^k,
\label{critC}
\end{equation}
where $\lambda_1^k$ denotes the first eigenvalue of \b{PCAmix} applied to the cluster $C_k$.

The first term (based on the squared Pearson correlation $r^2$) quantifies the link between  the numerical variables in $C_k$ and $\bbf^k$, independently of the sign of the relationship, whereas the second term (based  on the correlation ratio $\eta^2$) measures the link between the categorical variables in $C_k$  and $\bbf^k$. 
The homogeneity of a cluster is then maximized when all the numerical variables are perfectly correlated (or anti-correlated) to $\bbf^k$ and when all the correlation ratios of the categorical variables are equal to 1. 
In that case, all variables in cluster $C_k$ bring the same information which is summarized by the corresponding synthetic variable $\bbf^k$.

\medskip

\noindent
\textbf{Homogeneity ${\cal{H}}$ of a partition $P_K=\{C_1,\dots,C_K\}$}

\noindent
The homogeneity criterion ${\cal{H}}$ is defined as the sum of the homogeneities of its clusters:
\begin{equation}
{\cal{H}}(P_K)=\sum_{k=1}^K H(C_k)=\lambda_1^1+\ldots+\lambda_1^K.
\label{critP}
\end{equation}

\medskip

\noindent
\textbf{A hierarchical clustering algorithm}

\noindent
The objective of this algorithm is to find a partition of the $p$ available (numerical and/or categorical) variables. This partition must be  such that the variables within a cluster are strongly related to each other in the sense of the homogeneity criterion introduced previously. 
More specifically, for a given number $K$ of clusters, the aim is to find a partition $P_K$ which maximizes the homogeneity function ${\cal{H}}$ defined in (\ref{critP}). 
To this end, a hierarchical clustering algorithm can be used and is described hereafter.

This algorithm builds a set of $p$ nested partitions of variables as follows.
\begin{itemize}
\item {\sc Step $l=0$: Initialization.}  Start with the partition in $p$ clusters (i.e. one variable per cluster).
\item {\sc Step $l=1,\dots,p-2$: Aggregation of two clusters.} The objective is to aggregate two clusters of the partition  in $p-l+1$ clusters to get a new partition in $p-l$ clusters. 
To this end, we have to choose  the two clusters  $A$ and $B$ which provide the smallest dissimilarity $d(A,B)$ defined as:
\begin{eqnarray}
d(A,B)=H(A)+H(B)-H(A\cup B) =\lambda^A_1+\lambda^B_1-\lambda^{A\cup B}_1.  
\label{aggM}
\end{eqnarray}
This dissimilarity (aggregation measure) quantifies the lost of homogeneity observed when the two clusters $A$ and $B$ merge. 
Based on this measure of aggregation,  the new  partition in $p-l$ clusters maximizes ${\cal{H}}$ among all the partitions in $p-l$ clusters obtained by aggregation of two clusters of the partition in $p-l+1$ clusters.
\item {\sc Step $l=p-1$: Stop.} The final partition in one cluster (i.e. containing all the $p$ variables) is obtained.
\end{itemize}
%

\subsection{Variable selection using random forests}
\label{subsec:vsurf}

In this section, we focus on the \b{VSURF} procedure. This method is based
on random forests \citep{Breiman01}, which provide a non-parametric
predictor, with very good prediction performance in lots of applied
situations including high-dimensional data \citep[see e.g.][]{Verikas11}.
Implemented in the R package \texttt{VSURF} \citep{Genuer15}, the \b{VSURF}
procedure  is fully automatic (i.e. it does not
need a pre-specified number of variables to select) and can be applied for both
supervised classification and regression problems.
In this paper, we focus on the supervised classification case.

\medskip

\noindent
\textbf{Random Forests}

\noindent
A Random Forests (RF in the sequel) predictor is obtained by aggregating
a collection of randomly perturbed decision trees.
The randomness comes at two levels:
the individuals level with a preliminary bootstrap sample draw, and the variables level with random variables sub-samples draws before splitting a
node of a tree.


Let $\hat{h}_1, \ldots, \hat{h}_q$ be an ensemble of $q$ binary tree predictors.
The number of trees $q$ is usually set to several hundreds. Each tree is a piece-wise constant function, which associates a class label for every input vector.
A binary tree is made of internal nodes (which are split in two children nodes) and leaves (also called terminal nodes). The collection of the leaves forms a partition of the input space. The first node (which contains all data) is called the root of the tree, and a node is pure if it contains observations belonging to the same class.
Finally each internal node $t$ is split according to a splitting variable $X^{j_t}$ and a splitting value $v_t$. 
RF methodology is detailed in the following algorithm, which is implemented
in the {\tt randomForest} R package \citep{randomForest}.

{\sc Input:} a learning sample of size $n$ and an input vector $X$.

{\sc Goal:} predict the class label associated to $X$.

\begin{itemize}
  \item For each $l=1, \ldots, q$ :
  
    \begin{itemize}
      \item {\sc Bootstrap sample.} Draw a bootstrap sample of the learning set, by randomly choosing $n$ observations among the $n$ available, with replacement.
      \item {\sc Tree construction.}
      
        \begin{itemize}
          \item Initialize the tree by putting all observations of the bootstrap sample in the root node.
          \item {\it While it exists an impure node among current terminal nodes}, for each impure terminal node, randomly choose a subset of {\tt mtry} variables among the $p$ variables without replacement, seek the best split of the node only among splits involving the selected variables, and split the node into two children nodes according to the best split (see details below).
        \end{itemize}
      
      \item {\sc Tree prediction.}
      
      Let $X$ go down the tree and note the leaf it falls into. Return
      $\hat{h}_l(X)$, the class label of the observations of the learning set
      belonging to this leaf.
    \end{itemize}
    
  \item {\sc Aggregation.}
  
  Return the majority class among trees predictions:
  $\hat{f}(X) = \amax{c \in \{1, \ldots, L \} } \sum_{l = 1}^q \mathbf{1}_{\hat{h}_l(X) = c}$.
\end{itemize}

For supervised classification, the Gini index is used for heterogeneity measure of a node, in terms of class labels. It is defined, for a node $t$, as $\mathrm{Gini}(t) = \sum_{c \in \{1, \ldots, L \} } p_c(t) (1 - p_c(t))$, 
where $p_c(t)$ is the proportion of observations of class $c$ in node $t$.
%
Hence, the best split of a node $t$ is the one minimizing $\mathrm{Gini}(t)$ over all possible splits.
For a numerical variable, possible splits are of the form $\{ X^j \leq v \}$ with $v \in [ \min(\bx^j) , \max(\bx^j) [$. This split means that observations with a $j$-th variable value not larger than $v$ are sent to the left child node and the others to the right one.
For a categorical variable, possible splits are of the form $\{ X^j \in \M_j \}$ where $\M_j$ is a subset of the $j$-th variable levels (except empty set and total set).

Note that compared to Classification And Regression Trees (CART) introduced in \citet{Breiman84},
trees are fully developed
and are not pruned.

\medskip

\noindent
\textbf{OOB error and variable importance}

\noindent
During a RF run, an estimation of the prediction error and a measure of variable importance can be computed. Out-Of-Bag (OOB) error of a RF $\hat{f}$ is defined as $\mathrm{OOBerror}(\hat{f}) = \frac{1}{n} \mathrm{Card} \left\{ i \in \{ 1,\ldots,n \} \: | \: y_i \neq \hat{y_i} \right\}$,
where $\hat{y_i}$ is the majority class label among predictions of the trees $\hat{h}_l$ for which $y_i$ is OOB, that is for which $y_i$ was not chosen in the bootstrap sample used to build $\hat{h}_l$.
Variable Importance (VI) also uses the OOB samples (all observations of the learning set not included in a bootstrap sample) to compute a measure of the link between $Y$ and a variable $X^j$.
%
The idea is that the more the mean error of a tree on its OOB sample increases when the link between $X^j$ and $Y$ is broken, the more important the variable is.

\medskip

\noindent
\textbf{Variable Selection Using Random Forests}

\noindent
The \b{VSURF} procedure works with three steps. 
\begin{itemize}
  \item The first one begins by sorting variables based on random forests VI,
  and eliminates useless variables by an adaptive thresholding. The threshold
  is set to the estimation of the standard deviation (over multiple RF runs) of
  the VI of an unimportant variable.
  
  \item The second one starts with the
previously kept variables and performs an ascendant variable
introduction strategy, which builds embedded RF models. The
model which attains the minimum OOB error rate is then selected, and the
variables set on which it is based is called the interpretation
set.

 \item The third step consists in eliminating the redundancy of
interpretation variables and leads to a smaller variables set called
the prediction set. It consists in a step-wise ascendant strategy, which
at each step verifies that the next variable to introduce helps to decrease
enough the OOB error rate.
\end{itemize}
\noindent More details on those three steps can be found in \cite{Genuer15}.

\subsection{The CoV/VSURF procedure}
\label{subsec:covsurf}

We now describe the proposed methodology, which combines \b{CoV} and \b{VSURF},
in the following algorithm. First, the algorithm performs ``groups of
informative variables'' selection in (a). It encompasses the ascendant
hierarchical clustering of variables, the optimal choice of the number of clusters of
variables according to prediction performance, and the selection of the most important synthetic variables.
Then, in (b), the algorithm computes a prediction of a new observation $\bx$,
by building a classifier based on the groups of variables selected in (a).

\bigskip

{\sc Input:} a dataset ($\bX, \by$) and a new observation $\bx$ of $X$.

{\sc Goal:} select groups of informative variables and predict the class label
of $\bx$.

\begin{itemize}
\item[(a)] {\sc Groups of informative variables selection:} 

\begin{enumerate}
\item Apply \b{CoV} on $\bX$ to obtain a hierarchy (a tree) of variables.

\item For each $K = 2, \ldots, p$, cut \b{CoV} tree in $K$ clusters, train a RF
with the $K$ synthetic variables ${\bf f}^1, \ldots, {\bf f}^K$ as predictors
and $\by$ as output variable and compute its OOB error rate.

\item Choose the optimal number $K^*$ of clusters, which leads to the minimum
OOB error rate. 

Cut \b{CoV} tree in $K^*$ clusters.

\item Perform \b{VSURF} with the $K^*$ synthetic variables
${\bf f}^1, \ldots, {\bf f}^{K^*}$ as predictors and $\by$ as output variable.
Denote by $m \leq K^*$ the number of selected informative synthetic variables
(corresponding to the interpretation set of \b{VSURF}).
\end{enumerate}

\item[(b)] {\sc Prediction of a new observation $\bx$:}

\begin{enumerate}
  \item Train a random forest, $\hat{f}$, on the dataset consisting of the $m$
  selected synthetic variables and $\by$.

  \item Compute the scores of $\bx$ on the $m$ selected synthetic variables and
  predict its class label using $\hat{f}$.
\end{enumerate}

\end{itemize}

\bigskip

The main feature of this algorithm is that it outputs a list of $m$ selected
informative synthetic variables. Since each synthetic variable is built on a subset
of original variables, the algorithm implicitly leads to select groups of original variables, the group structure being a priori unknown. Hence,
in addition to perform variable selection, it takes advantage of the
clustering of variables to give a list of selected variables with the
additionnal group structure information.

Another interesting feature of the procedure is that even if the ascendant
hierarchical clustering of variables algorithm is \emph{unsupervised} (in the
sense that it does not use response variable $\by$), the final variables
partition is \emph{supervised} since the number of clusters is optimized
in terms of RF classifier prediction error. This choice is justified by
the fact that our main goal is prediction, so we are primarily interested
in groups of informative variables, rather than groups of all variables
(e.g., a complex group structure only involving non-informative variables does
not have to be found by the method to give good prediction).
We stress that each group of informative variables is summarized by its 
synthetic variable (the first principal component of the group).
These synthetic variables are then used to build the predictive model.
Typically, informative variables should be well represented by their associated
synthetic variables. Indeed, assume that the information of a predictive
variable $X^j$ is not captured by the first principal component of its group
$C$ (for instance $X^j$ is orthogonal---not correlated---to the component).
The procedure would then select a higher number of clusters ($K^*$) in order
to split $C$ into two new informative groups $C_1$ and $C_2$, with e.g. 
$C_1$ containing variables highly correlated to the synthetic variable
associated to $C$, and $C_2$ with a synthetic variable that retrieves the
prediction power of $X^j$. 

This procedure is illustrated on simulated samples and a real proteomic
dataset, respectively  in Section~\ref{sec:simu-study} and Section~\ref{sec:real-data}. Finally, let us mention that the method
\b{CoV/VSURF} has been implemented in an R package, which is available online
together with a vignette guide
\footnote{\url{https://github.com/robingenuer/CoVVSURF/blob/master/vignettes/intro_CoVVSURF.Rmd}}.

\section{Simulation study}
\label{sec:simu-study}

In this section, we evaluate numerical performance of the \b{Cov/VSURF}
procedure in the framework of supervised classification.
In the following simulation model, explanatory variables are structured in
informative or non-informative groups. The underlying groups of variables
can be numerical, categorical or mixed. Let us recall that \b{Cov/VSURF}
does not take into account a priori the group structure information.
We first describe the simulation model used to generate the data. We then consider an ``$n > p$" (resp. an ``$n < p$") simulated dataset to evaluate the proposed methodology in comparison with alternative approaches.

\subsection{Simulation model}

Let us consider a binary response variable $Y$ and a $p$-dimensional explanatory variable $X$.
The conditional probability $P(Y=1|X=x)=:p(x)$ is modeled as a function of $x$ using the well known logistic regression model:
\begin{equation}\label{model-sim}
\log\left(\frac{p(x)}{1-p(x)}\right)=x'\beta-9,
\end{equation}
where $\beta\in \R^p$.
Note that the term ``-9'' in (\ref{model-sim}) is used to center the index $X'\beta$ (for the considered choice of parameters below) and therefore to obtain equibalanced response for $Y$.

Let us now specify how we construct the $p$-dimensional explanatory variable $X$ (with $p=120$)  in order to get $10$ groups of variables (components $X^j$ of $X$) of several types (numerical, categorical or mixed) such that  some  groups of variables are informative or not. Table~\ref{tabsim} provides a complete description of the 10 groups of these $p=120$ variables.
To sum up, we have a binary response variable and nine groups of correlated variables: one small and one large informative groups of numerical variables, one moderate non-informative group of numerical variables, and this structure is repeated for the three groups of categorical variables and the three groups of mixed variables. In addition, we have $30$ non-informative and non-structured (independent) variables in the last group. Details on data generating process are given below.

\begin{table}[!ht]
\caption{\label{tabsim} Overview of the $10$ groups of variables in the simulation study.}
\centering
\begin{tabular}{ccccc}
  \toprule
 Type  & Size  & Informative & Components  & Names   \\
 of variables & of group & group &  of $\beta$ &  of variables  \\
  \midrule
 Numerical &Small (3) &  Yes & (1, 2, 3)& {\tt NumS$j$} with $j=1,\dots,3$ \\
 Numerical &Moderate (12)& No & $(0,\dots,0)$ &  {\tt NumM$j$} with $j=1,\dots,12$ \\
 Numerical &Large (15) & Yes & $\frac{1}{5}(\underbrace{1,
  \ldots, 1}_{5}, \underbrace{2, \ldots, 2}_{5},
\underbrace{3, \ldots, 3}_{5})$& {\tt NumL$j$} with $j=1,\dots,15$ \\
  Categorical &Small (3)& Yes & (1, 2, 3)& {\tt CategS$j$} with $j=1,\dots,3$  \\
Categorical&Moderate (12) & No &$(0,\dots,0)$&{\tt CategM$j$} with $j=1,\dots,12$ \\
Categorical&Large (15) & Yes &$\frac{1}{5}(\underbrace{1,
  \ldots, 1}_{5}, \underbrace{2, \ldots, 2}_{5},
\underbrace{3, \ldots, 3}_{5})$& {\tt CategL$j$} with $j=1,\dots,15$ \\
 Mixed&Small (3) & Yes &(1, 2, 3) &{\tt MixedS$j$} with $j=1,\dots,3$ \\
 Mixed&Moderate (12) & No &$(0,\dots,0)$&{\tt MixedM$j$} with $j=1,\dots,12$ \\
 Mixed&Large (15) & Yes &$\frac{1}{5}(\underbrace{1,
  \ldots, 1}_{5}, \underbrace{2, \ldots, 2}_{5},
\underbrace{3, \ldots, 3}_{5})$& {\tt MixedL$j$} with $j=1,\dots,15$\\
Numerical &Big (30) & No &$(0,\dots,0)$& {\tt Noise$j$}  with $j=1,\dots,30$ \\
 \bottomrule
 \end{tabular}
\end{table}

Let $Z$ be the following multivariate Gaussian distribution ${\mathcal{N}}_p (\mu, \Sigma)$
where $p=120$, $\mu = 0_p$ and the covariance matrix $\Sigma$ is a block-diagonal matrix in order to get $10$ groups of variables. These groups are independent from each other.
Let us also introduce the parameter $\rho \in [-1, 1]$ to control the link between the variables within a group.
 The covariance matrix of $Z$ is then defined as:
\begin{equation}
\Sigma = \mathrm{diag}\left( \Sigma_{3, \rho}, \Sigma_{15, \rho},
  \Sigma_{12, \rho}, \Sigma_{3, \rho}, \Sigma_{15, \rho},
  \Sigma_{12, \rho}, \Sigma_{3, \rho}, \Sigma_{15, \rho},
  \Sigma_{12, \rho}, \sigma^2 I_{30} \right),
\label{sigma-def}
\end{equation}
where
$\Sigma_{s,  \rho}$ is an $s\times s$ matrix, whose coefficients are all equal to $\rho$ except on the diagonal where they are equal to $1$, $\sigma^2>0 $ and $I_{30}$ is the identity matrix of dimension $30$.
When $\rho$ is close to $1$, we obtain nine groups of highly correlated numerical variables and one group of $30$ variables which are independent. In the following, we set $\rho=0.9$.

%
%
We then generate $n$ mutually independent random $p$-dimensional vectors $\bz_1,\dots,\bz_n$ from the Gaussian distribution ${\mathcal{N}}_p(\mu, \Sigma)$.
From the $\bz_i$'s, we first construct
 the output variable values $y_i \in \{0, 1\}$ using the logistic regression model (\ref{model-sim}) and $\beta$ coordinates given in Table~\ref{tabsim}.
Since the  $\bz_i$'s original data are numerical, we have to binarize some numerical variables $Z^j$ to obtain categorical ones by
thresholding at their median value: if the variable value is less than
the median, the value of the corresponding binarized variable $X^j$  is $0$, and this value is $1$ otherwise.
In the 3 groups of mixed variables, one third of the variables are binarized (the last, last 4, last 5 variables in the small, medium and large group respectively).


\subsection{An ``$n > p$'' simulated dataset}
\label{sec:large-simul-datas}

\noindent
\textbf{Application of \b{CoV/VSURF}  to a learning dataset} 

\begin{figure}[!ht]
\begin{center}
\includegraphics[width=1.2\textwidth,angle=270]{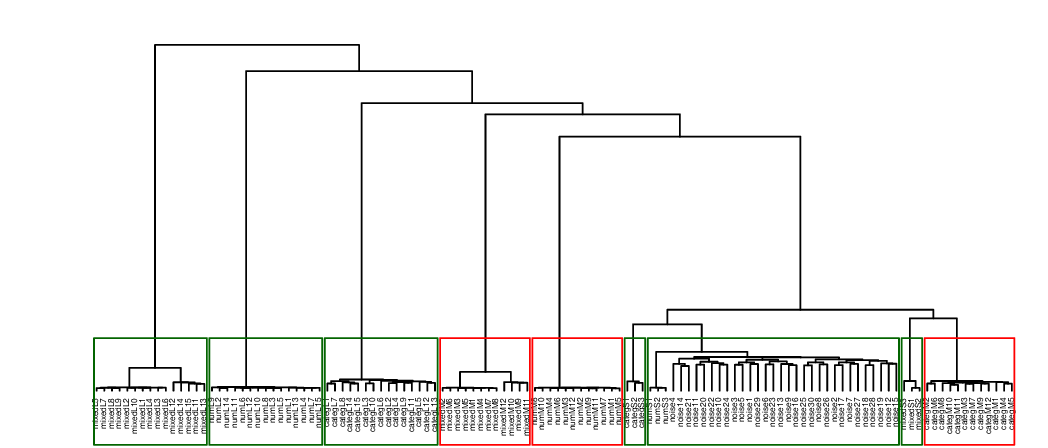} 
\caption{\b{CoV} tree of the $p=120$ variables based on the learning sample of size $n=600$. Heights correspond to dissimilarities between two aggregated clusters. Rectangles indicate the partition in $K^*=9$ clusters (in green, the $m = 6$ clusters selected by \b{VSURF}, and in red the remaining unrelevant clusters).}
\label{fig:dendro_600}
\end{center}
\end{figure}

\begin{figure}[!ht]
\begin{center}
\includegraphics[width=0.45\textwidth]{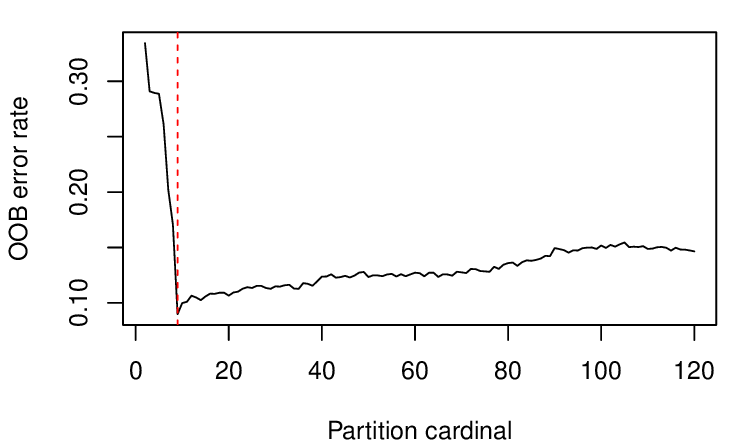} 
\end{center}
\caption{Random Forests OOB error rate according to the number of clusters of the partitions obtained by cutting the \b{CoV} tree for the simulated learning dataset of $n=600$ observations. The dashed red vertical line corresponds to $K^* = 9$ clusters.}\label{fig:<res_covsurf_don600}
\end{figure}

\noindent
From the previous model, a learning sample of $n=600$ observations of the $p=120$ explanatory variables and of the binary response variable $Y$ is generated. 
The method \b{CoV}  is applied to the $600 \times 120$ matrix of explanatory variables. The  dendrogram  of the hierarchy of the 120 variables  (also called the \b{CoV} tree)  is given in Figure \ref{fig:dendro_600}. This dendrogram suggests a partition in $9$ clusters.
However, in our methodology, the number of clusters  is not chosen according to the shape of the dendrogram but according to the prediction of the binary response variable $Y$. For each value of $K$ between $2$ and $120$, we cut the \b{CoV} tree, build a random forest on the $K$ synthetic variables of the corresponding $K$-clusters partition and compute the RF OOB error rate. This procedure is illustrated Figure \ref{fig:<res_covsurf_don600}. The (nearly) optimal value we get for this learning dataset is $K^*=9$. This partition in $K^*=9$ clusters
retrieves almost the complete structure of the data. We recover 8 of the 9 groups of  correlated variables, while all noise variables are pulled
together in a large cluster with the last group of correlated variables (the small one).  However, the synthetic variable of this cluster does not really take the noise variables into account, whose loadings (coefficients in the linear combination) are very low contrary to the loadings of the variables of the small group of correlated variables. Since this small group of numerical variables is informative, it means that this synthetic variable is also informative.

 \b{VSURF} is then applied on the $K^*=9$ synthetic variables of the
previously chosen partition. \b{VSURF} selects 6 synthetic variables,
corresponding to the $6$ informative groups of variables.
From the interpretation point of view, we
succeed in selecting all informative variables, with in addition
the clustering structure.
For comparison, \b{VSURF} directly applied on the 120 original variables selects 39 variables among the 54 informative ones, with at least one per informative group (the large group of categorical variables being the least recovered).

\medskip

\noindent
\textbf{Prediction performances on a test dataset}

\begin{figure}[!ht]
\begin{center}

\begin{minipage}{0.45\textwidth}
\includegraphics[width=\textwidth]{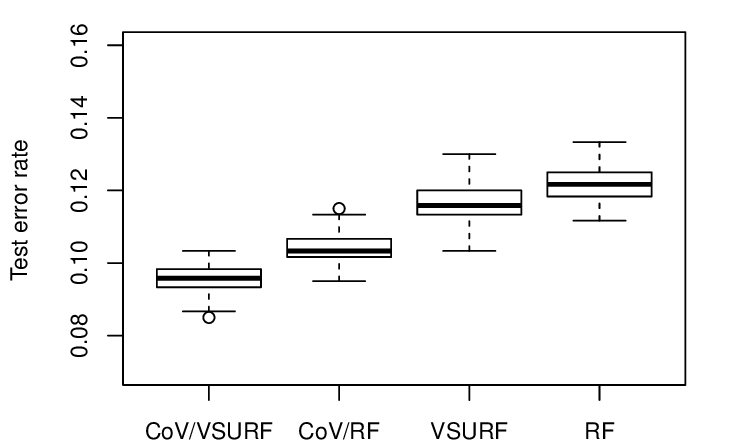} \\
\centerline{(a) Single learning dataset}
\end{minipage}
\begin{minipage}{0.45\textwidth}
\includegraphics[width=\textwidth]{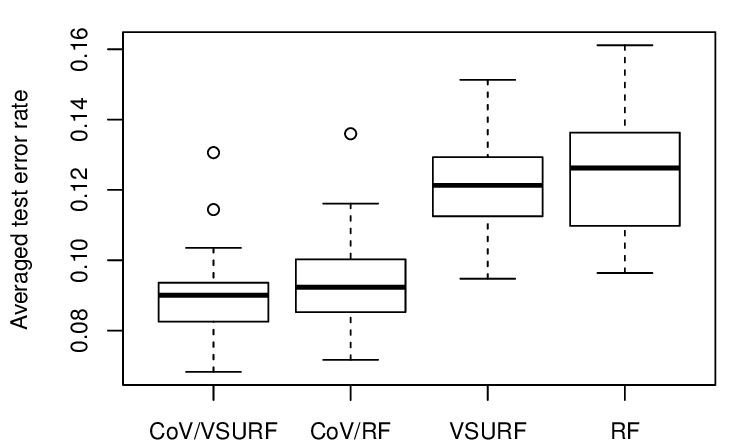} \\
\centerline{(b) Replication of 50 learning datasets}
\end{minipage}

\end{center}

\caption{Comparison of test error rates with learning datasets of size $n = 600$.
\b{CoV/VSURF} and \b{CoV/RF} correspond to \b{VSURF} and random forests (RF) applied on the $K^*$ synthetic variables given by \b{CoV}. \b{VSURF} and \b{RF} refer to methods applied on the original 120 variables. (a) Test error rates of 100 forests trained on one learning dataset. (b) Averaged (over 100 forests) test error rates computed on 50 learning datasets.}
\label{fig:boxplot_don600}
\end{figure}

\noindent
In addition to the learning dataset, a test sample of $n=600$ observations is now generated.
We focus on the prediction performances of the proposed methodology \b{CoV/VSURF}  on this test sample and do comparisons with 3 other methodologies: 
\begin{itemize}
\item \b{VSURF}:  variable selection using random forests  is applied on the original 120 variables,
\item  \b{RF}: random forests are applied  on the original 120 variables,
 \item \b{CoV/RF}: random forests are applied on the $K^*$ synthetic variables obtained by \b{CoV}.
 \end{itemize}
Note that the test sample of size $n=600$ is fixed all along this simulation study.
We obtain in Figure~\ref{fig:boxplot_don600} the boxplots of test error rates.

On the left-hand side of Figure~\ref{fig:boxplot_don600},  boxplots correspond to 100 runs of random forests. More precisely,
\begin{itemize}
\item  For \b{CoV/VSURF} and \b{CoV/RF}, the clustering of variables is performed once and the optimal number $K^*$ of clusters is chosen automatically using OOB error rate. For \b{Cov/RF}, 100 forests are trained on the $K^*$ synthetic variables, while for \b{Cov/VSURF}, 100 forests are trained on the $m < K^*$ synthetic variables selected by \c{VSURF} (which is then performed also a single time). 
\item For \b{VSURF}, 100 forests are trained on the variables selected with one run of \c{VSURF} applied on the $p=120$ original variables. 
\item For \b{RF}, 100 forests are trained on the original variables.
\end{itemize}
The results are better with  \b{CoV/VSURF} and  \b{CoV/RF} i.e. when \b{VSURF} and \b{RF} are applied after the dimension reduction using \b{CoV}. In addition, the approach \b{CoV/VSURF} is slightly better than  \b{CoV/RF}  in this setting. Also, the variability due to the random nature of the \b{RF} is relatively small here and very similar with the four approaches.

On the right-hand side of Figure~\ref{fig:boxplot_don600}, 50 learning datasets are generated in order to take into account the variability of the results due to the random nature of the simulation procedure.
The four methodologies explained above for a single learning dataset are applied on each dataset. Boxplots correspond then to the 50 averaged (over 100 forests) test error rates.
The results are better here with  \b{CoV/VSURF} and  \b{CoV/RF}.  The
improvement obtained by grouping variables together before the classifier
run is confirmed here, despite the variability due to the random simulation
procedure.

\medskip

To sum up, using \b{CoV} is in a dimension reduction step leads  --- at least for this example ---  to a gain in prediction both with
\b{RF}  and \b{VSURF}. Moreover applying \b{VSURF} in \b{CoV/VSURF} procedure permits to select
informative groups of linked variables without loss in prediction.
In the next section,   the exact same experiment is repeated but in a more
challenging situation where only $n=60$ ($<p=120$) observations are generated.

\subsection{An ``$n < p$'' simulated dataset}
\label{sec:small-simul-datas}

\begin{figure}[!ht]
\begin{center}
\includegraphics[width=1.2\textwidth,angle=270]{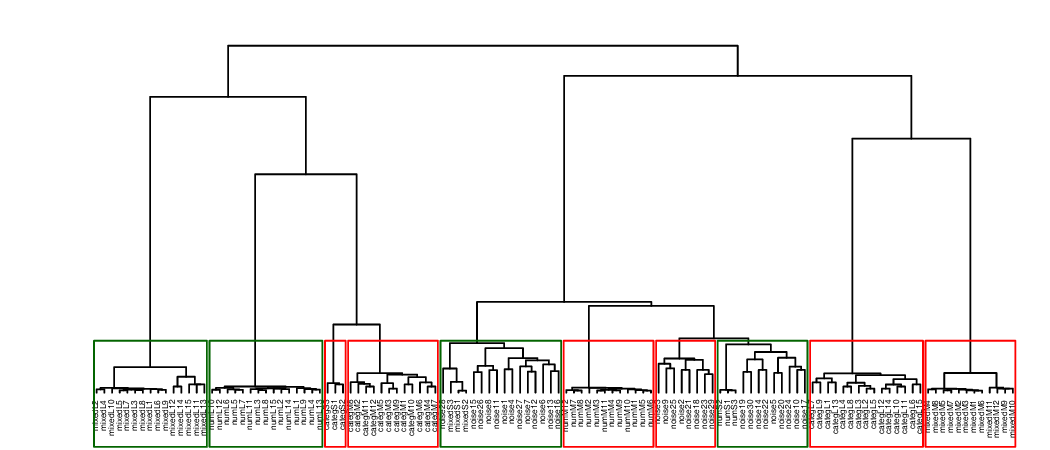} 
\caption{\b{CoV} tree of the $p=120$ variables based on the learning sample of size $n=60$. Heights correspond to dissimilarities between two aggregated clusters. Rectangles indicate the partition in $K^*=10$ clusters (in green, the $m = 4$ clusters selected by \b{VSURF}, and in red the remaining unrelevant clusters).}
\label{fig:dendro_60}
\end{center}
\end{figure}

Now, we simulate a learning sample of $n=60$ observations of the $p=120$ explanatory variables and of the binary response variable $Y$. 
The \b{CoV} method  is applied to the $60 \times 120$ matrix of  explanatory variables. Note that we set the number of trees $q$ to 2000 in all our experiments.
The dendrogram of the hierarchy is represented in Figure \ref{fig:dendro_60}.
This tree suggests to retain the partition in $K=9$ clusters which matches with the underlying partition in 9 groups of correlated variables (i.e. the informative and non informative structured groups of variables) while the 30 noise variables are pulled together either with the small group of numerical variables or with the small group of mixed variables.

\begin{figure}[!ht]
\begin{center}
\includegraphics[width=0.45\textwidth]{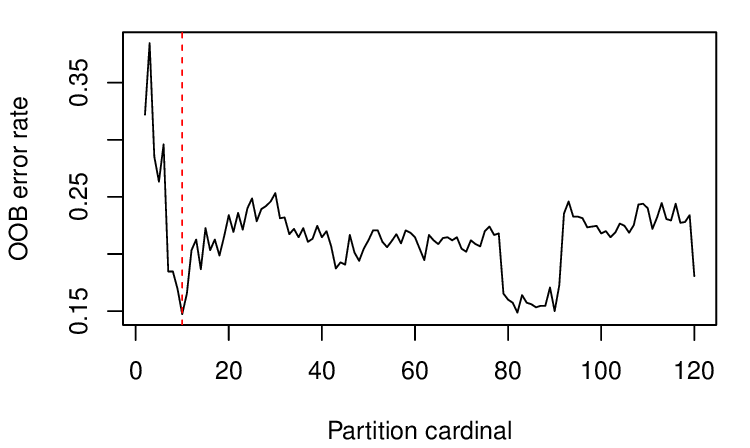} 
\end{center}
\caption{Random Forests OOB error rate according to \b{CoV} partition cardinal for the learning dataset of $n=60$ observations. The dashed red vertical line corresponds to $K^* = 10$ clusters.}\label{fig:res_covsurf_don60}
\end{figure}

Our methodology chooses to cut the dendrogram in $K^*=10$ clusters (see Figure~\ref{fig:res_covsurf_don60}). The structure in 9 groups of correlated variables is conserved and the 10-th cluster only contains 8 noise variables.
\b{VSURF} applied on the $K^*=10$ associated synthetic variables selects $m=4$ of them. The two groups of informative numerical variables and the two groups of informative mixed variables are retained.

Thus, contrary to the case with $n = 600$ observations (see Figure~\ref{fig:dendro_600}), where the group structure is perfectly identified and all informative variables selected,
in this more difficult context where $n < p$, the \b{CoV/VSURF} methodology still retrieves the group structure in the dendrogram (see Figure~\ref{fig:dendro_60}) but misses the informative groups of categorical variables in the selection step. This can be explained by the fact that those variables are binary, hence less explanatory than the numerical ones.
\b{VSURF} applied on the 120 original variables only selects 8 variables, with at least one per group selected by our approach, and still no informative categorical variable are retained.

\begin{figure}[!ht]
\begin{center}
\begin{minipage}{0.45\textwidth}
\includegraphics[width=\textwidth]{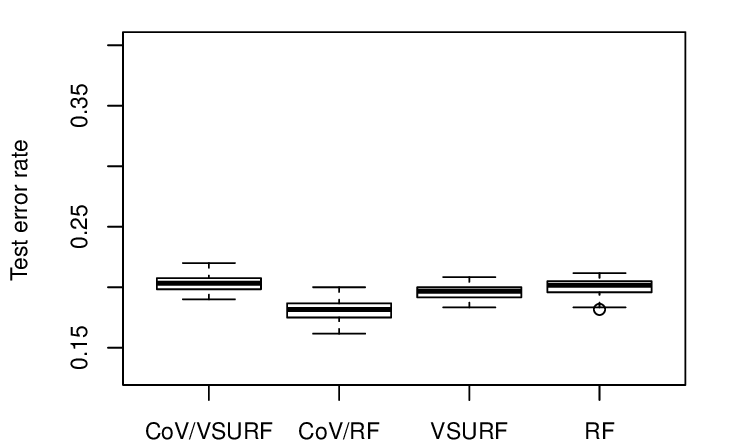} \\
\centerline{(a) Single learning dataset}
\end{minipage}
\begin{minipage}{0.45\textwidth}
\includegraphics[width=\textwidth]{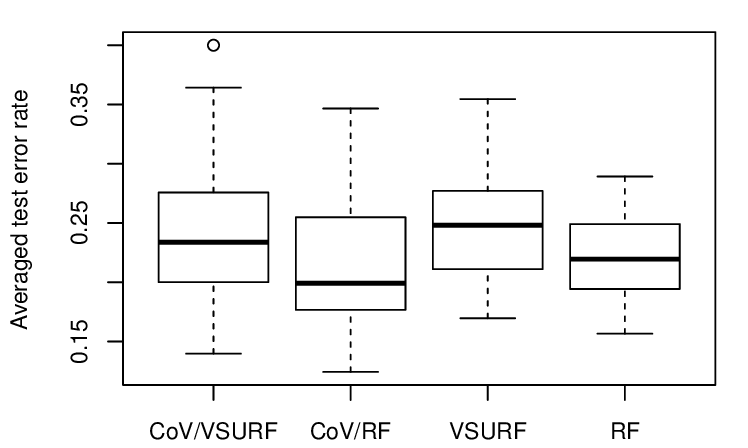} \\
\centerline{(b) Replication of 50 learning datasets}
\end{minipage}
\end{center}
\caption{Comparison of test error rates with learning datasets of size $n = 60$.
\b{CoV/VSURF} and \b{CoV/RF} correspond to \b{VSURF} and random forests (RF) applied on the $K^*$ synthetic variables given by \b{CoV}. \b{VSURF} and \b{RF} refer to methods applied on the original 120 variables. (a) Test error rates of 100 forests trained on one learning dataset. (b) Averaged (over 100 forests) test error rates computed on 50 learning datasets.}
\label{fig:boxplot_don60}
\end{figure}

To study prediction performance of the proposed methodology \b{CoV/VSURF}, the  same test sample as in the previous section is used, and Figure~\ref{fig:boxplot_don60} is obtained in the same way than in the $n=600$ case.
Those graphs show that the variability due to random forests (left-hand side) remains relatively small compared to the variability due to random generation of the 50 learning dataset (right-hand side).
In this more difficult situation, all methodologies are fairly comparable in prediction, with a slight advantage to \b{CoV/RF}. However, the methodology \b{CoV/VSURF} has the advantage to select groups of informative variables, without increasing prediction error rate.

\subsection{Sensitivity analysis}
\label{sec:sensitivity-analysis}

To conclude this simulation study, we conduct a sensitivity analysis with
respect to the simulated group structure. Indeed, in the previous simulation
model, the covariance matrix $\Sigma$ (see Equation~\ref{sigma-def}) induces
a very clear group structure in the data, which could be somehow unrealistic
for real applications. Thus, we relax this structure by adding a
between-groups correlation coefficent $\delta \geq 0$. More precisely, the
intra-groups correlaction structure is unchanged, but two variables from two
different groups have now a correlation of $\delta$, while noise variables
remain independent with every other variables.

\begin{figure}[!ht]
\begin{center}
\begin{minipage}{0.45\textwidth}
\includegraphics[width=\textwidth]{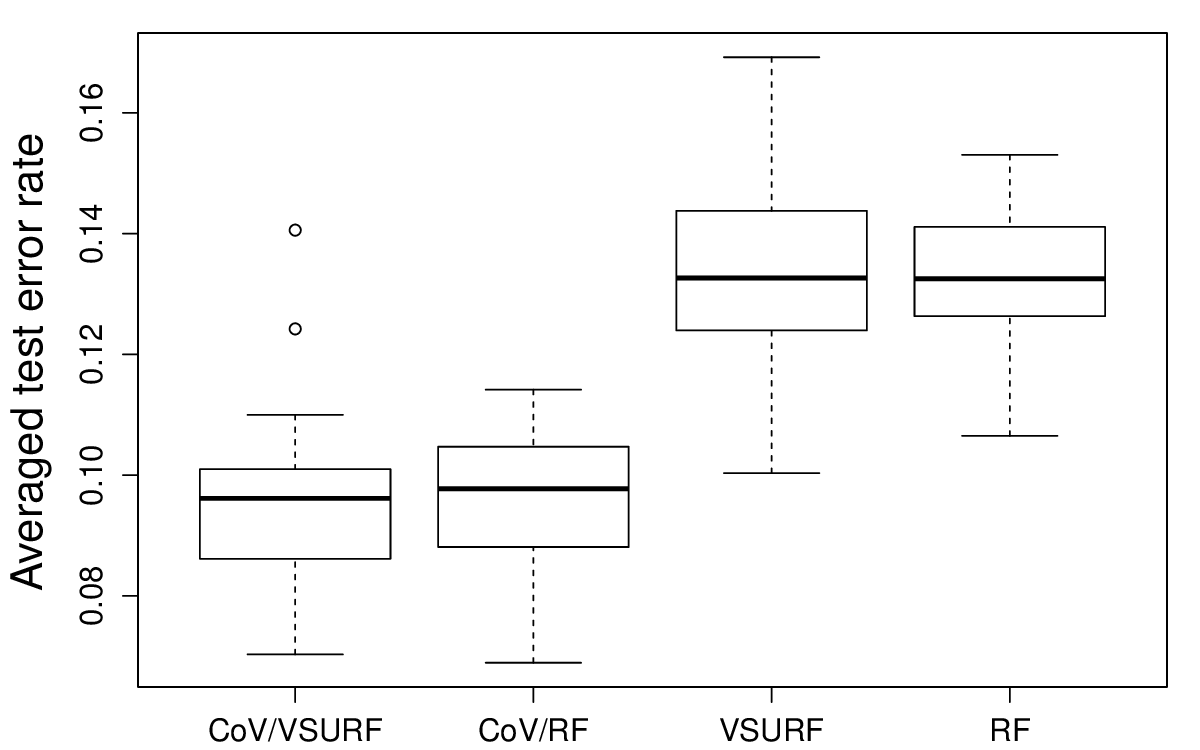} \\
\centerline{(a) $n = 600, \delta = 0.3$}
\end{minipage}
\begin{minipage}{0.45\textwidth}
\includegraphics[width=\textwidth]{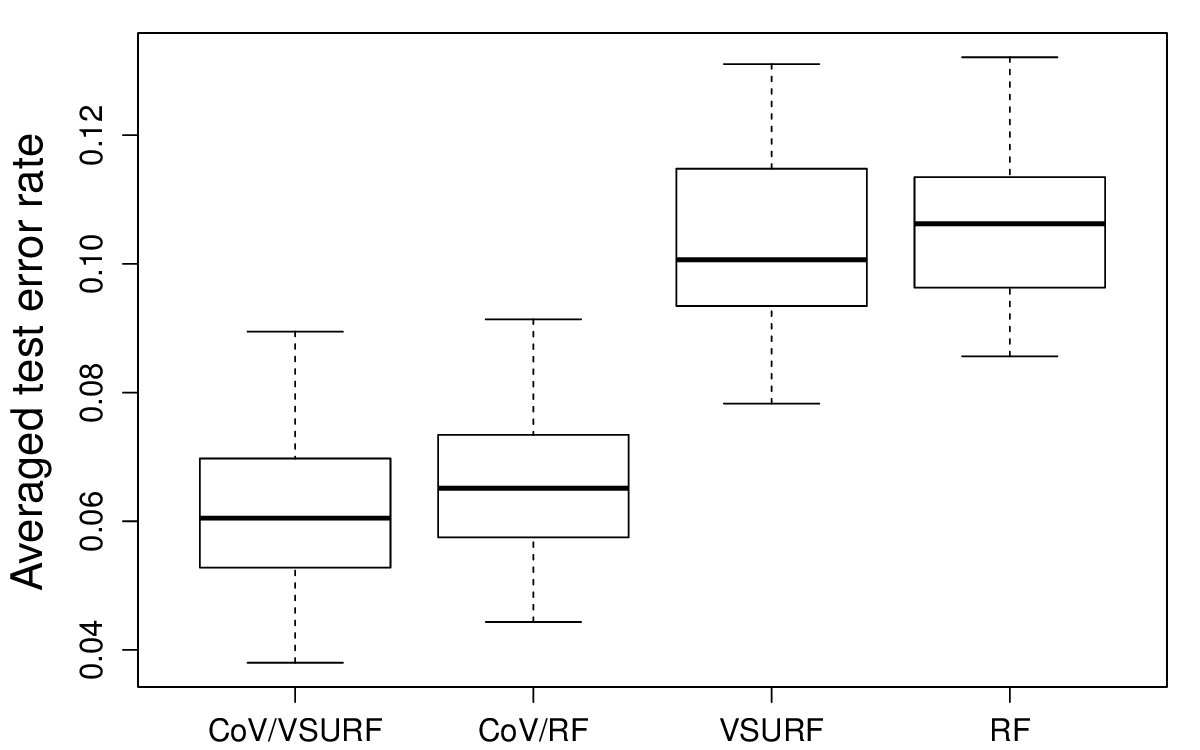} \\
\centerline{(b) $n = 600, \delta = 0.6$}
\end{minipage}
\begin{minipage}{0.45\textwidth}
\includegraphics[width=\textwidth]{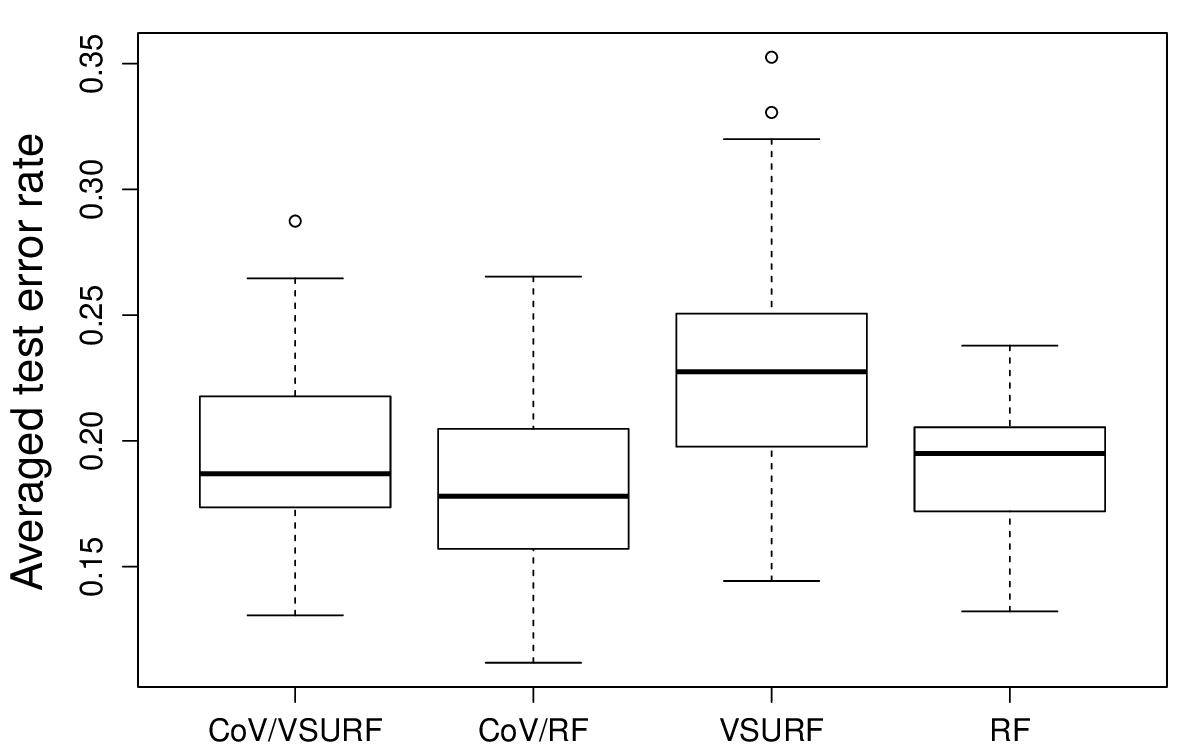} \\
\centerline{(c) $n = 60, \delta = 0.3$}
\end{minipage}
\begin{minipage}{0.45\textwidth}
\includegraphics[width=\textwidth]{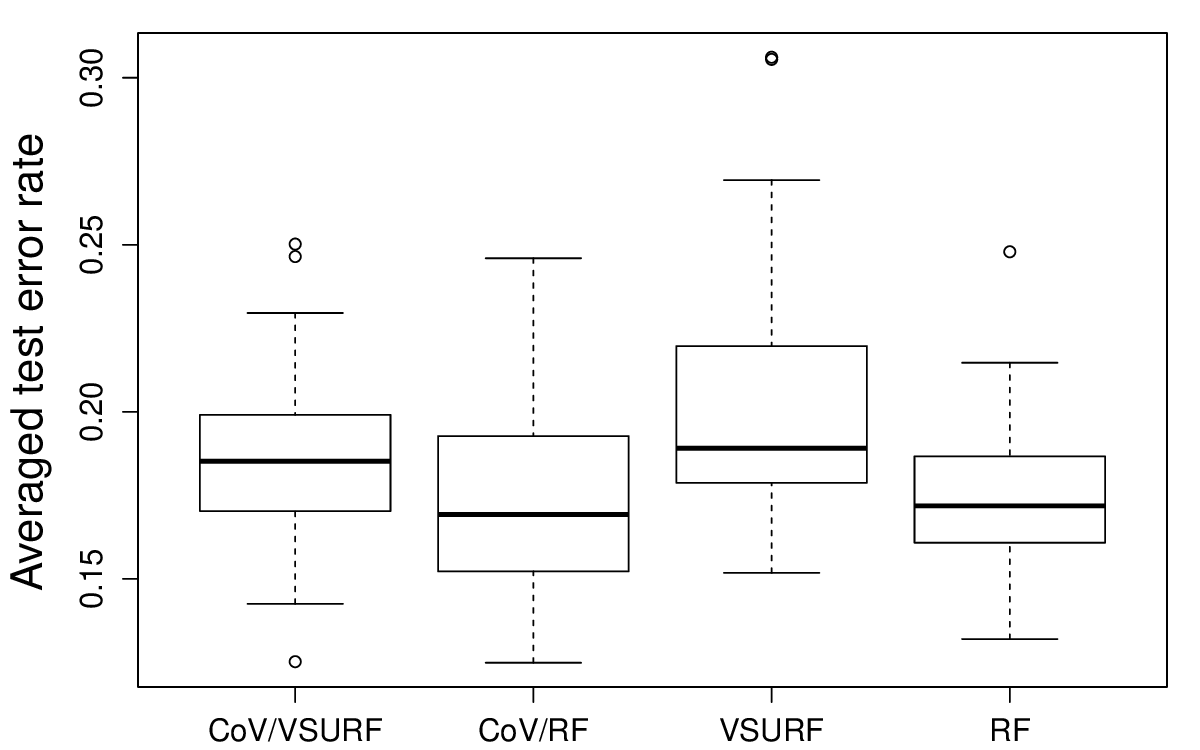} \\
\centerline{(d) $n = 60, \delta = 0.6$}
\end{minipage}
\end{center}
\caption{Comparison of test error rates with learning datasets of size $n = 600$
(top) and $n = 60$ (bottom) including a between-groups correlation of
$\delta = 0.3$ (left) and $\delta = 0.6$ (right). Boxplots correspond
to the averaged (over 100 forests) test error rates computed on 50 learning
datasets.}
\label{fig:sensitivity}
\end{figure}

Figure~\ref{fig:sensitivity} illustrates the results in the two previsous cases
with $n=600$ and $n=60$, where $\delta$ takes the values $0.3$ and $0.6$. Boxplots correspond to the 50 averaged (over 100 forests)
test error rates associated to 50 randomly generated datasets as in
Figures~\ref{fig:boxplot_don600}(b) and \ref{fig:boxplot_don60}(b).

Those results show that the comparison between the four methods
\b{CoV/VSURF}, \b{CoV/RF}, \b{VSURF} and \b{RF} is quite stable, with
a clear advantage in terms of prediction for methods using clustering of
variables in the case $n=600$ and similar performances of the four methods
in the case $n=60$.

Furthermore, we investigate the results of the \b{CoV/VSURF} procedure applied
to one dataset (results not shown here) for several values of $\delta$. It
reveals that the method is not perturbed by values of $\delta$ up
to 0.6 (the clustering of variables is accurate and the selection of synthetic variables is satisfactory), and that for larger values it naturally begins to give very small number
of groups. For example in the extreme case of $\delta = 0.9$, the method groups
together all variables (except noise variables) as expected.

\section{Proteomic data application}
\label{sec:real-data}

In this section,   the \b{CoV/VSURF} procedure is illustrated with an application
on a real proteomic dataset. This dataset
comes from a clinical study sponsored by University Hospital in Bordeaux.
It involved $n=44$ patients with a rectum cancer who
undertook a treatment of chemotherapy and radiotherapy, before a
surgery intervention. Some patients responded favorably to the treatment
and hence had a smaller tumor at $t_1$, the time of surgery, compared to $t_0$,
a time just before the beginning of the treatment, and some patients did not.
The main goal of this study was to predict if the patient will be a good
 treatment responder or not, using proteomic information, measured at $t_0$.
Indeed, it is useless to give the treatment to a bad responder, for which other alternatives should be tried. A secondary objective is to select the proteins which best discriminate the two kinds of treatment response.

The dataset contains $p=4786$ numerical explanatory variables which measure protein abundances. Measurements are done on peptides (one protein is made of several peptides). The total number of peptides is $p=4786$, whereas total number of proteins is $868$. A priori knowledge of what peptide is part of which protein is not used in our approach, but it will be used for interpretation purpose.
This helps us to interpretate the resulting group structure
and permits to illustrate the interest of our method for other problems
where no a priori group structure is known.

\begin{figure}[!ht]
\begin{center}
 \includegraphics[width=0.45\textwidth]{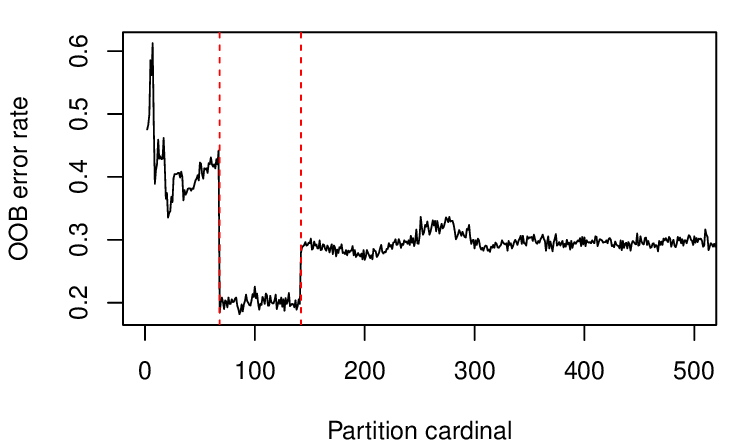} 
\caption{Random forests OOB error rate according to \b{CoV} partition cardinal (number of peptides clusters) for the proteomic data. The two dashed red vertical lines correspond respectively to partitions in $K^*=68$ and $K=142$ clusters. For visibility, the x-axis is truncated around $K=500$ clusters (the error rate remains quite stable until $K=2000$).} \label{fig:dendro_proteo}
\end{center}
\end{figure}

The \b{CoV/VSURF} methodology is then applied to this dataset.
Since the number $p$ of variables is relatively large, the number $K^*$ of clusters of peptides is chosen among 2 and 2000. The resulting random forests OOB error rates are displayed in Figure~\ref{fig:dendro_proteo} and the (nearly) optimal number of clusters is $K^*=68$.
The \b{VSURF} method, applied on those $K^*=68$ synthetic variables, selects 4 of them.
These 4 synthetic variables are sorted by decreasing order of their variable importance (VI).
The corresponding 4 groups of peptides gather respectively 37, 61, 20, 25 peptides,
that is a total of 143 peptides. Those 143 peptides come from
73 different proteins. The most homogeneous one is the 4th with 22
peptides (over 25) coming from the same protein, \textit{PZP} (Pregnancy Zone Protein).
The 3rd one contains 18 peptides (over 20) coming from two different proteins, \textit{IGHD}  (Ig Delta Chain) protein and \textit{THRB} (Prothrombin) protein.
The second cluster is the more heterogeneous with 61 peptides coming from 50 different proteins.
The first one contains 13 peptides (over 37) coming from \textit{A2GL} (Alpha 2 Glycoprotein), 6 from \textit{CRP} (C-reactive protein).

The cluster with highest VI, made of 37 peptides,
is the one responsible of the large decrease of OOB error rates in
Figure~\ref{fig:dendro_proteo}: it appears when the number of
synthetic variables goes from $67$ to $68$ i.e. when the OOB error rate decreases
from around $40\%$ to $20\%$. In addition, when this particular cluster
is split into two new clusters (when the partition cardinal $K$ reaches 142),
the OOB error increases significantly from $20\%$ to $30\%$.
This confirms the importance of this cluster of peptides to predict treatment
response.
Note that, among the peptides selected by \b{CoV/VSURF},
we find peptides that come from \textit{APOA1} (Apolipoprotein A1), \textit{HPT}
(Haptoglobin), \textit{TRFE} (Serotransferrin) and \textit{PGRP2}
(N-Acetylmuramoyl L-Alanine Amidase). Those proteins make sense to oncologists
in this context, who are currently investigating in more depth the
relation between those proteins and the response to the treatment.
In addition, when \b{VSURF} is applied alone on the original data,  only 35 peptides are selected. So, \b{VSURF} gives a sparse variable selection (less peptides are selected), but with no group structure.

\begin{figure}[!ht]
\begin{center}
\includegraphics[width=0.5\textwidth]{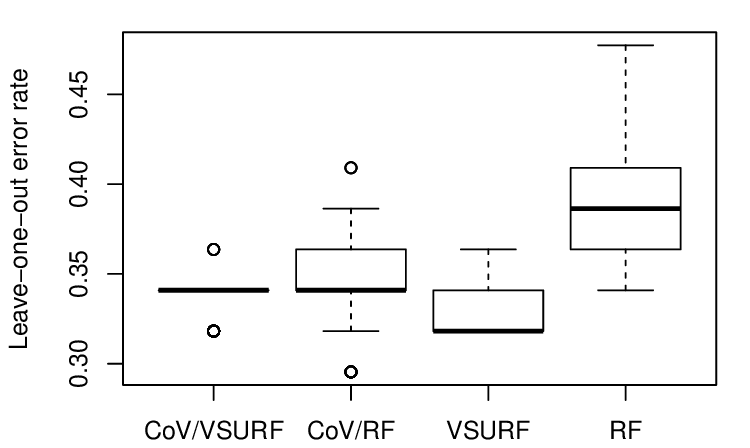} 
\caption{Comparison of leave-one out error rates of 100 forests for proteomic data. \b{CoV/VSURF} and \b{CoV/RF} correspond to \b{VSURF} and random forests (RF) applied on the $K^*=68$ synthetic variables given by \b{CoV}. \b{VSURF} and \b{RF} refer to methods applied on the original $p=4786$ peptides.}\label{fig:boxplot_proteo}
\end{center}
\end{figure}

To get an estimation of classification error rates for this real
dataset, we perform an external leave-one-out cross-validation
to fairly estimate prediction performance of our approach \citep[see][]{Ambroise02}.
This means that, for each patient in this study, we perform the entire methodology \b{CoV/VSURF} on the dataset containing all other patients, before predicting him.
As before, we compare the performance with the three approaches \b{CoV/RF}, \b{VSURF} and \b{RF}.
Figure~\ref{fig:boxplot_proteo} provides the results.
According to these boxplots, there is almost no difference between the performances of  \b{CoV/VSURF}, \b{CoV/RF} and \b{VSURF}
(recall that since $n = 44$, one more good classification of an individual gives approximately a $0.02$ decrease in error rate).
It is important to highlight that  only the \b{CoV/VSURF} procedure selects groups of informative peptides without increasing  the prediction error. To conclude this real data study, \b{CoV/VSURF} achieves the double objective of predicting  if the patient is a good treatment responder and  determining groups of peptides  responsible for the good or the bad response to the treatment.

\section*{Acknowledgments.} 


The authors would like to thank members of the USMR unit of Bordeaux University Hospital (especially H\'{e}l\`{e}ne Savel and Paul Perez) for allowing access to the clinical study data and for insightful discussions.
In addition, the authors also thank the editor and two anonymous reviewers
who contribute to greatly improve the manuscript.

\bibliographystyle{apalike}
\bibliography{chavent_genuer_saracco.Cov_VSURF_VersionArxivV3}

\end{document}